\documentclass{article}%
\usepackage{amssymb}
\usepackage{amsfonts}%
\usepackage{amsmath}%
\usepackage{graphicx}
\newcommand {\bd} {\begin{displaymath}}
\newcommand {\ed} {\end{displaymath}}
\newcommand {\be} {\begin{equation}}
\newcommand {\ee} {\end{equation}}

\newtheorem{theorem}{Theorem}

\begin{document}

\title{  A Theorem  involving the denominators of Bernoulli numbers}

\author{  Pantelis A.~Damianou, Peter Schumer}
\date{}
\maketitle

\begin{abstract}
Consider the average of the first $n$ $k$-th powers. We pose and answer the following natural question: For which values of $n$ and $k$ is this average an integer? If $k$ is odd the answer is easy; it is an integer as long as $n$ is incongruent to $2$ modulo $4$. If $k$ is even then the criterion involves the denominator of the $k$-th Bernoulli number. The average is an integer iff $n$ is not divisible by any prime which divides the denominator of the $k$-th Bernoulli number.
\end{abstract}

The Swiss mathematician, Jakob Bernoulli (1654 - 1705),
successfully sought a general method for summing the first n k$^{th}$ powers
for arbitrary positive integers n and k. \ Let us define

\bd  S_{k}(n) = \sum_{j=1}^{n}j^{k} = 1^{k} + 2^{k} +  \cdots  +
n^{k} \ . \ed

Define the average of the first $n$  $k$th powers by

\bd  \mu_k(n)=\frac{S_k(n) }{n}  \ . \ed

We pose and and answer the following natural question: For which values of $n$ and $k$ is $\mu_k(n)$ an integer?. Our answer, although it does involve the denominators of Bernoulli numbers, which undergraduates may not have see, relies primarily upon elementary divisibility arguments.
\bigskip

\noindent
{\bf Background} In his \textit{Ars Conjectandi}, published posthumously in 1713 and dedicated
primarily to the theory of probability, Bernoulli  presented a recursive solution for
S$_{k}$(n). \ It states that for k $\geq$ 1,

\bd
(n + 1)^{k+1} = (n + 1) + \sum_{j=1}^{k}\binom{k+1}{j}S_{j}(n)   \ , \ed

where the binomial coefficients are defined as usual
 \bd \binom{k+1}{j} = \frac{(k+1)!}{j!(k+1-j)!}  \ . \ed

Furthermore, if we define what are now called the \textit{Bernoulli numbers} by

\bigskip

\bd  B_{0} = 1,  \qquad   (k + 1)B_{k} = - \sum_{j=0}^{k-1}\binom{k+1}{j}B_{j}   \ \ \ \ \ \ \ \ j \ge 1   \ed

then for $k \ge 1$, the sums $S_{k}(n)$  satisfy:

\bd
(k + 1) S_{k}(n) = \sum_{j=0}^{k}\binom{k+1}{j} B_{j} (n +
1)^{k+1-j}  \ .  \ed

The Bernoulli numbers are the rational coefficients of the linear terms of the
(k+1)$^{st}$ degree polynomials, $S_{k}(n-1)$. For example,

\bd  S_{0}(n -1) = 1n- 1 \ed
\bd  S_{1} (n -1) = \frac{1}{2}n^{2} - \frac{1}{2}n  \ed
\bd  S_{2} (n -1) =\frac{1}{3}n^{3} - \frac{1}{2}n^{2} + \frac{1}{6}n \ed
\bd  S_{3}(n -1) =\frac{1}{4}n^{4} - \frac{1}{2}n^{3} + \frac{1}{4}n^{2} + 0n \ed

\bd S_{4} (n -1) = \frac{1}{5}n^{5} - \frac{1}{2}n^{4} - \frac{4}{15}
n^{3} -\frac{1}{30}n \ .  \ed

\bigskip

 It follows that \ B$_{0}$ = 1, B$_{1}$ =
-$\frac{1}{2}$, B$_{2}$ = $\frac{1}{6}$, B$_{3}$ = 0, and B$_{4}$ = -$\frac
{1}{30}$. \ In fact, B$_{2n+1}$ = 0 for all n $\geq$ 1. \ More compactly, we
can define the Bernoulli numbers by the following power series:

\bd
\frac{x}{e^{x}\text{ - 1}} = \sum_{k=0}^{\infty}\frac{B_{k}\text{ x}^{k}
}{\text{m!}} \ . \ed

\bigskip
For even $k\ge 2$ , we write B$_{k}$ = N$_{k}$/D$_{k}$, where N$_{k}$ and
D$_{k}$ are relatively prime and D$_{k}$ $\geq$ 1. The numerators N$_{k}$ have
played a significant role in number theory due largely to their connection
with Fermat's Last Theorem. A prime $p$ is a \textit{regular} prime if $p$ does
not divide any of the numbers N$_{2}$, N$_{4}$, ..., N$_{p-3}$. (The only
irregular primes less than 100 are 37, 59, and 67.) In 1850, Ernst Kummer
proved that Fermat's Last Theorem is true for every exponent which is a
regular prime. Of course, Andrew Wiles (1995) has since proven Fermat's Last
Theorem in toto.

The denominators D$_{k}$ have played a less significant role in mathematics
even though they can be clearly described. The Von Staudt - Clausen Theorem
(1840) states that for even $k$,  D$_{k}$ is the product of all primes $p$ with (p -1)%
%TCIMACRO{\TEXTsymbol{\vert}}%
%BeginExpansion
$\vert$%
%EndExpansion
k. An interesting consequence is that D$_{k}$ is square-free for all $k$.  The
theorem was proven independently (and nearly simultaneously) by the two mathematicians.

\bigskip
\bigskip

\textbf{Examples  }

We begin by considering a few examples, deriving results  directly using congruence relations.

$\bigskip$

\begin{itemize}
\item   k = 1: We have $\mu_{1}(n) = \frac{n+1}{2}$. Hence $\mu_{1}(n)  \in \mathbb{Z} $  iff $n$  is odd. This is an exceptional case due to the fact that $B_1 \not=0$.

\item  $k = 2$: In this case, $\mu_{2}$(n) = $\frac{(n+1)(2n+1)}{6}$. We claim
that $\mu_2(n) \in \mathbb{Z} $  iff $n$  is not divisible by $2$ or $3$. First, suppose that $n$ is not divisible by
$2$ or $3$. Clearly, $(n +1)(2n +1)$ is even. If  $n \equiv 1 \ ({\rm mod}\, 3)$,  then 3%
%TCIMACRO{\TEXTsymbol{\vert}}%
%BeginExpansion
$\vert$%
%EndExpansion
(2n +1) and if  $n \equiv 2 \ ({\rm mod}\, 3)$, then 3%
%TCIMACRO{\TEXTsymbol{\vert}}%
%BeginExpansion
$\vert$%
%EndExpansion
(n +1). In any event, 6%
%TCIMACRO{\TEXTsymbol{\vert}}%
%BeginExpansion
$\vert$%
%EndExpansion
(n +1)(2n +1) and so $\mu_{2}(n) \in \mathbb{Z}$. Second, suppose that $n$ is divisible by either $2$ or $3$. If 2%
%TCIMACRO{\TEXTsymbol{\vert}}%
%BeginExpansion
$\vert$%
%EndExpansion
n, then (n +1)(2n +1) is odd and hence is not divisible by 6. If 3%
%TCIMACRO{\TEXTsymbol{\vert}}%
%BeginExpansion
$\vert$%
%EndExpansion
n, then $n = 3k$ for appropriate integer $k$, and $(n +1)(2n +1) = (3k +1)(6k +1) =
18k^{2} + 9k + 1$, a number not divisible by $3$ (nor by $6$).

\item $ k = 3$: We have $\mu_{3}$(n) = $\frac{n(n+1)^{2}}{4}$. We claim that
$\mu_{3}(n) \in \mathbb{Z} $  as long as $n$ is incongruent to $2$ modulo $4$. If $n$ is congruent to $0$, $1$, or $3$
modulo $4$, then 4%
%TCIMACRO{\TEXTsymbol{\vert}}%
%BeginExpansion
$\vert$%
%EndExpansion
n(n +1)$^{2}$. However, if n $\equiv$ $2$ \ (mod $4)$, then n(n +1)$^{2}$
$\equiv$ $2$(mod $4)$, and so 4 does not divide $n(n+1)^{2}$.

\item $ k = 4$: In this case, $\mu_{4}$(n) = $\frac{(n+1)(2n+1)(3n^{2}%
+3n-1)}{30}$. We claim that $\mu_{4}(n) \in \mathbb{Z}$  iff $n$ is not divisible by $2$, $3$, or $5$. Suppose that $n$ is relatively prime to
$30$  (equivalently not divisible by $2$, $3$, or $5$). Then $n +1$ is even and $(n +1)(2n
+1)$ is divisible by $3$. Furthermore, (n +1)(2n +1)(3n$^{2}$ + 3n -1) = 6n$^{4}$
+ 15n$^{3}$ + 10n$^{2}$ - 1 $\equiv$ n$^{4}$ - 1 (mod 5). But by Fermat's
Little Theorem, n$^{4}$ -1 $\equiv$ 0 (mod 5) and so 5%
%TCIMACRO{\TEXTsymbol{\vert}}%
%BeginExpansion
$\vert$%
%EndExpansion
(n +1)(2n +1)(3n$^{2}$ + 3n -1). Hence $\mu_{4}(n) \in \mathbb{Z}$   in this case.

In the other direction, if 2%
%TCIMACRO{\TEXTsymbol{\vert}}%
%BeginExpansion
$\vert$%
%EndExpansion
n, then (n +1)(2n +1)(3n$^{2}$ + 3n -1) is odd and not divisible by 30. If 3%
%TCIMACRO{\TEXTsymbol{\vert}}%
%BeginExpansion
$\vert$%
%EndExpansion
n, then (n +1)(2n +1)(3n$^{2}$ + 3n -1) $\equiv$ -1 (mod 3) and so is not
divisible by 30. Finally, if 5%
%TCIMACRO{\TEXTsymbol{\vert}}%
%BeginExpansion
$\vert$%
%EndExpansion
n, then (n +1)(2n +1)(3n$^{2}$ + 3n -1) $\equiv$ -1 (mod 5) and so is not
divisible by 30.

\end{itemize}

These examples hint that the situation is very different for odd and even values of $n$.  We develop our main theorem in two sections. Only the even case involves the Bernoulli numbers. In both parts, we use the easily noted fact that $\mu_k(n) $ is an integer if and only if $S_k(n) \equiv  0 \ ({\rm mod} \  n) $.

$\bigskip$

\noindent
\textbf{An odd Theorem:}

\begin{theorem}  For odd numbers $k\ge 3$,    $\mu_{k}$(n) is an integer iff n is incongruent to 2
modulo 4.
\end{theorem}

{\it Proof} \ Suppose $k$ is odd and $k \ge 3$.  Since (n - a)$^{k}$ $\equiv$ -a$^{k}$ (mod n) for all a, we can pair up the
terms of $S_{k}(n)$ from the outside in.

$\qquad$(a) If $n$ is odd, then

\bd
S_{k}(n) = [1^{k} + (n -1)^{k}] + [2^{k} + (n - 2)^{k}] + ... +
[\left(\frac{n-1}{2}\right)^{k} + \left(\frac{n+1}{2}\right)^{k}] + n^{k} \equiv  \ed

\bd
= (1^{k} - 1^{k}) + (2^{k} - 2^{k}) + ... + 0 \equiv 0 \  ({\rm mod} \,  n) \ .
\ed

\qquad(b) If n is even, then there are two subcases depending on whether or
not n is divisible by 4.

\qquad\qquad(i) If n $\equiv$ 0 (mod 4), then

\bd
S_{k}(n) = [1^{k} + (n -1)^{k}] + [2^{k} + (n - 2)^{k}] + ... +
[\left(\frac{n}{2} -1\right)^{k} + \left(\frac{n}{2} +1\right)^{k}] + (\frac{n}{2})^{k}
+ n^{k}
\equiv 0 \ \ ({\rm mod} \  n)  \ed
since $k >1$ and $\frac{n}{2}$ is even.
\bigskip

\qquad\qquad(ii) If n $\equiv$ 2 (mod 4), then

\bd S_{k}(n) = [1^{k} + (n -1)^{k}] + [2^{k} + (n - 2)^{k}] + ... +
[\left(\frac{n}{2} -1\right)^{k} + \left(\frac{n}{2} +1\right)^{k}] + (\frac{n}{2})^{k}
+ n^{k} \ed

\bd
\equiv (\frac{n}{2})^{k} \ \  ({\rm mod} \  n). \ed

But $\frac{n}{2}$ is odd and so ($\frac{n}{2}$)$^{k}$ is odd. Since n is even,
($\frac{n}{2}$)$^{k}$ is incongruent to 0 (mod n).

\bigskip
\noindent
{\bf An even more interesting theorem}

\begin{theorem}
For even numbers $k \ge 2$,  $\mu_{k}$(n) is an integer iff $n$ is relatively prime to $D_k$.

\end{theorem}

\noindent
{\it Proof} \   The Von Staudt-Clausen Theorem ([\textbf{1}], Theorem 118) states that the k$^{th}$
Bernoulli denominator D$_{k}$ = $\Pi_{(p-1)|2k}$ p (see [\textbf{1}],
 Since in this case k is even, we may rewrite the formula as D$_{k}$ =
$\Pi_{(p-1)|k}$ p. To prove our result it must be shown that S$_{k}$(n)
$\equiv$ 0 (mod n) iff for every prime p%
%TCIMACRO{\TEXTsymbol{\vert}}%
%BeginExpansion
$\vert$%
%EndExpansion
n that p$\nmid\ $D$_{k}$. By Von Staudt-Clausen it suffices to establish that

\begin{center}
S$_{k}$(n) $\equiv$ 0 (mod n) iff for every prime p%
%TCIMACRO{\TEXTsymbol{\vert}}%
%BeginExpansion
$\vert$%
%EndExpansion
n that (p -1)$\nmid$k.\qquad\qquad(1)
\end{center}

We will utilize the following easily established result valid for any prime p
([\textbf{1}], Theorem 119):

\begin{center}
$\sum_{m=1}^{p}$ m$^{k}$ $\qquad\equiv$ -1 (mod p) if (p -1)%
%TCIMACRO{\TEXTsymbol{\vert}}%
%BeginExpansion
$\vert$%
%EndExpansion
k

\qquad\qquad$\qquad\qquad\qquad\qquad\qquad\equiv$ 0 (mod p) if (p -1)$\nmid
$k\qquad\qquad\qquad(2)

\bigskip
\end{center}

It is convenient to first prove the theorem  assuming \textit{n square-free}. \

We establish (1):

($\Longleftarrow$) Suppose that for all p%
%TCIMACRO{\TEXTsymbol{\vert}}%
%BeginExpansion
$\vert$%
%EndExpansion
n that (p -1)$\nmid$k. Choose a prime p%
%TCIMACRO{\TEXTsymbol{\vert}}%
%BeginExpansion
$\vert$%
%EndExpansion
n. By (2),

\begin{center}
$\sum_{m=1}^{p}$ m$^{k}$ $\equiv$ 0 (mod p).
\end{center}

Similarly,

\begin{center}
$\sum_{m=rp+1}^{(r+1)p}$ m$^{k}$ $\equiv$ 0 (mod p) for 0 $\leq$ r $\leq$
$\frac{n}{p}$ -1.
\end{center}

Hence S$_{k}$(n) = $\sum_{m=1}^{n}$ m$^{k}$ = $\sum_{r=0}^{\frac{n}{p}-1}%
\sum_{m=rp+1}^{(r+1)p}$ m$^{k}$ $\equiv$ 0 (mod p) .

But p arbitrary and n square-free implies that S$_{k}$(n) $\equiv$ 0 (mod n).

($\Rightarrow$) Suppose there exists a prime p%
%TCIMACRO{\TEXTsymbol{\vert}}%
%BeginExpansion
$\vert$%
%EndExpansion
n such that (p -1)%
%TCIMACRO{\TEXTsymbol{\vert}}%
%BeginExpansion
$\vert$%
%EndExpansion
k. By (2)

\begin{center}
$\sum_{m=1}^{p}$ m$^{k}$ $\equiv$ -1 (mod p).
\end{center}

Similarly,

\begin{center}
$\sum_{m=rp+1}^{(r+1)p}$ m$^{k}$ $\equiv$ -1 (mod p) for 0 $\leq$ r $\leq$
$\frac{n}{p}$ -1.
\end{center}

Hence S$_{k}$(n) $\equiv$ -$\frac{n}{p}$ (mod p), which is incongruent to 0
(mod p) since p and $\frac{n}{p}$ are relatively prime. Thus S$_{k}$(n) is
incongruent to 0 (mod n).

\bigskip

Now suppose that \textit{n is not square-free}.

($\Longleftarrow$) Suppose that for all p%
%TCIMACRO{\TEXTsymbol{\vert}}%
%BeginExpansion
$\vert$%
%EndExpansion
n that (p -1)$\nmid$k. If there is a prime $p$ exactly dividing $n$ (that is $p|n$, but $p^2$ does not divide $n$), Then as in the square-free case, $S_k(n) \equiv 0 \ \ ({\rm mod} \, p)$.  Now let $p$ be specifically a prime $p$ with p$^{a}$%
%TCIMACRO{\TEXTsymbol{\vert}}%
%BeginExpansion
$\vert$%
%EndExpansion%
%TCIMACRO{\TEXTsymbol{\vert}}%
%BeginExpansion
$\vert$%
%EndExpansion
n with $a \geq 2$. (The notation p$^{a}$%
%TCIMACRO{\TEXTsymbol{\vert}}%
%BeginExpansion
$\vert$%
%EndExpansion%
%TCIMACRO{\TEXTsymbol{\vert}}%
%BeginExpansion
$\vert$%
%EndExpansion
n means that p$^{a}$%
%TCIMACRO{\TEXTsymbol{\vert}}%
%BeginExpansion
$\vert$%
%EndExpansion
n and p$^{a+1}\nmid$n.)

\bigskip

Lemma: Let p be a prime with (p -1)$\nmid$k. Then

\begin{center}
1$^{k}$ + 2$^{k}$ + ... + (p$^{a}$)$^{k}$ $\equiv$ 0 (mod p$^{a}$).

\bigskip
\end{center}

Proof of Lemma: (Induction on a)

If a = 1, then

\begin{center}
1$^{k}$ + 2$^{k}$ + ... + p$^{k}$ $\equiv$ 0 (mod p) by (2).
\end{center}

Assume then that the lemma holds for a - 1, namely that

\begin{center}
1$^{k}$ + 2$^{k}$ + ... + (p$^{a-1}$)$^{k}$ $\equiv$ 0 (mod p$^{a-1}$).
\end{center}

Now consider S$_{k}$(p$^{a}$) = $\sum_{r=0}^{p-1}\sum_{j=1}^{p^{a-1}}$
(rp$^{a-1}$ + j)$^{k}$.

The binomial theorem implies that

\begin{center}
(rp$^{a-1}$ + j)$^{k}$ = $\sum_{i=0}^{k}$ $\binom{k}{i}$r$^{i}$ p$^{(a-1)i}$
j$^{k-i}$.
\end{center}

Hence

\begin{center}
S$_{k}$(p$^{a}$) = $\sum_{r=0}^{p-1}\sum_{j=1}^{p^{a-1}}$ $\sum_{i=0}^{k}$
$\binom{k}{i}$r$^{i}$ p$^{(a-1)i}$ j$^{k-i}$.\qquad(3)
\end{center}

For i $\geq$ 2, p$^{(a-1)i}$ $\equiv0$ (mod p$^{a}$) and so all terms of (3)
with i $\geq$ 2 are congruent to 0 (mod p$^{a}$).

For i = 0, $\sum_{r=0}^{p-1}\sum_{j=1}^{p^{a-1}}$ j$^{k}$ = p$\cdot$S$_{k}%
$(p$^{a-1}$).

But S$_{k}$(p$^{a-1}$) $\equiv$ 0 (mod p$^{a-1}$) by our inductive hypothesis.
Hence $\sum_{r=0}^{p-1}\sum_{j=1}^{p^{a-1}}$ j$^{k}$ $\equiv$ 0 (mod p$^{a}$).

For i = 1,

\begin{center}
$\sum_{r=0}^{p-1}\sum_{j=1}^{p^{a-1}}$ krp$^{a-1}$j$^{k-1}$ = $\sum
_{r=0}^{p-1}$ krp$^{a-1}\cdot$S$_{k-1}$(p$^{a-1}$)

= k S$_{k-1}$(p$^{a-1}$)$\cdot$p$^{a-1}\cdot\frac{(p-1)p}{2}$.
\end{center}

But $S_{k-1}(p^{a-1}) \in \mathbb{Z} $  and 2%
%TCIMACRO{\TEXTsymbol{\vert}}%
%BeginExpansion
$\vert$%
%EndExpansion
(p -1). Thus

\begin{center}
$\sum_{r=0}^{p-1}\sum_{j=1}^{p^{a-1}}$ krp$^{a-1}$j$^{k-1}$ $\equiv$ 0 (mod
p$^{a}$).
\end{center}

Therefore, S$_{k}$(p$^{a}$) $\equiv$ 0 (mod p$^{a}$) and the lemma is proven.

\bigskip

Analogous to the lemma, it follows that

\begin{center}
$\sum_{m=rp^{a}+1}^{(r+1)p^{a}}$ m$^{k}$ $\equiv$ 0 (mod p$^{a}$) for 0 $\leq$
r $\leq$ $\frac{n}{p^{a}}$ - 1.
\end{center}

Hence S$_{k}$(n) $\equiv$ 0 (mod p$^{a}$). But since p was arbitrary, S$_{k}
$(n) $\equiv$ 0 (mod n).

\bigskip

($\Longrightarrow$) A slight modification of the square-free proof works here, as follows.
On the one hand, if there  exists a prime p%
%TCIMACRO{\TEXTsymbol{\vert}}%
%BeginExpansion
$\vert$%
%EndExpansion
n such that (p -1)%
%TCIMACRO{\TEXTsymbol{\vert}}%
%BeginExpansion
$\vert$%
%EndExpansion
k and p$^{a}$%
%TCIMACRO{\TEXTsymbol{\vert}}%
%BeginExpansion
$\vert$%
%EndExpansion%
%TCIMACRO{\TEXTsymbol{\vert}}%
%BeginExpansion
$\vert$%
%EndExpansion
n with a $\geq$ 2. By (2), $\sum_{m=1}^{p}$ m$^{k}$ $\equiv$ -1 (mod p). Hence
S$_{k}$(n) $\equiv$ -$\frac{n}{p}$ which is incongruent to 0 (mod p$^{a}$)
since p$^{a}$%
%TCIMACRO{\TEXTsymbol{\vert}}%
%BeginExpansion
$\vert$%
%EndExpansion%
%TCIMACRO{\TEXTsymbol{\vert}}%
%BeginExpansion
$\vert$%
%EndExpansion
n. Thus n$\nmid$S$_{k}$(n) as in the square-free case .

On the other hand, suppose there exists a prime $p$ with  p$^{a}$%
%TCIMACRO{\TEXTsymbol{\vert}}%
%BeginExpansion
$\vert$%
%EndExpansion%
%TCIMACRO{\TEXTsymbol{\vert}}%
%BeginExpansion
$\vert$%
%EndExpansion
n
with $a \ge 2$ and (p -1)%
%TCIMACRO{\TEXTsymbol{\vert}}%
%BeginExpansion
$\vert$%
%EndExpansion
k. By (2), $\sum_{m=1}^{p}$ m$^{k}$ $\equiv$ -1 (mod p). Thus
\bd \sum_{m=1}^p m^k \equiv (rp-1) \ ({\rm mod}\, p^a) \  \ed
for some $r$ with  $1 \le r \le p^{a-1}$. But then $S_k(n) \equiv \frac{n}{p} (rp-1) \equiv -\frac{n}{p}  \ ({\rm mod} p^a) $. Hence $S_k(n)$ is not congruent to 0  modulo $p^a$ and so n$\nmid$S$_{k}$(n).  This completes the proof of the theorem.
 $\square$

\begin{center}
\bigskip
\end{center}

Reference:

1. G.H. Hardy and E.M. Wright, \textit{An Introduction to the Theory of
Numbers }(5th ed.), Oxford, Clarendon Press, 1979.

\bigskip
\noindent\textit{Department of Mathematics and Statistics,
University of Cyprus, P.O.~Box 20537, 1678,  Nicosia, Cyprus\\
damianou@ucy.ac.cy}

\bigskip
\noindent\textit{Department of Mathematics, Middlebury College,
Warner Hall 306, Middlebury, VT 05753, \\
schumer@middlebury.edu}
\end{document}